\documentclass[review]{elsarticle}

\usepackage{lineno,hyperref}
\usepackage{natbib}
\setcitestyle{authoryear,open={((},close={))}}
\modulolinenumbers[5]
\usepackage{amsmath}
\usepackage{bbm}
\usepackage{amsthm}
\usepackage{amssymb}
\usepackage{antpolt}
\usepackage[polish]{babel}
\usepackage{polski}
\newtheorem{defn}{Definition}
\newtheorem{theorem}{Theorem}
\newtheorem{lemma}{Lemma}
\begin{document}
	\begin{frontmatter}
		
		\title{On some Limit Theorem for Markov Chain}
		
		\author{Anna Czapkiewicz}
		\address{	Faculty of Management,\\
			AGH University of Science and Technology\\
	gzrembie@cyf-kr.edu.pl}
		
		\author{Antoni Leon Dawidowicz}

		\address{Faculty of Mathematics and Computer Science,\\	Jagiellonian University}
	\begin{abstract}
	 The goal of this  paper is to describe conditions which
	 guarantee a central limit theorem for random variables, which distributions are controled by hidden 
	 Markov chains. We proved that when  a Markov chain is ergodic and random variables fullfiled Lindeberg's condition  then the Central Limit Theorem is true.
	\end{abstract}
	
	
\end{frontmatter}

\section{Introduction}
Regime switching models have been used extensively in econometric time series analysis. In most of these models, two regimes are introduced with a state process determining one of the regimes to take place in each period. The bivalued state process is typically modeled as a Markov chain. The autoregressive model with this type of Markov switching was first considered by Hamilton (1989), and later analyzed by Kim (1994). 
  Markov-switching models with endogenous explanatory variables have
 been considered  by Kim (2004, 2009).  The most authors  assume that the Markov chain, which determine regimes, is completely independent from all other parts of the model. Diebold et al. (1994) and Kim (2009)   considers a Markov-switching driven by a set of observed variables. 
Chang et al. (2016) introduces a new approach to model regime switching using an autoregressive latent factor, which determines regimes depending upon whether it
takes a value above or below some threshold level. 

Despite numerous generalizations of this type of models, there is still little known about their theoretical properties. For example, one of the problems  is the likelihood ratio test and other tests for comparing two regime switching models, the second is the indication of regularity conditions for the  efficiency of maximum likelihood estimator of unknown model parameters. 
Various statistical properties of the model have been studied by Hansen (1992), Hamilton (1996), Garcia (1998), Timmermann (2000), and Cho and White (2007), among others. The  overview of the  literature is in  monograph by Kim and Nelson (1999).
In order to solve many problems related to testing hypothesis or some estimator efficiency, it is enough to prove Central Limit Theorem.
In regime switching models the random variables, which distributions are controled by hidden 	 Markov chains we have dependent variables. 

The central limit theorem has been extended to the case of
dependent random variables by several authors. The conditions under which these theorems are stated either are very restrictive or involve conditional distributions, which makes them difficult to apply. Hoeffding and Robbins (1994) prove central limit theorems for sequences of dependent random variables of a certain special type which occurs frequently in mathematical statistics. 

In this paper we prove Central Limit Theorem for random variables, which distributions are controled by hidden 	 Markov chains. We prove that when  a Markov chain is ergodic and random variables fullfiled Lindeberg's condition  then the Central Limit Theorem is true.
\section{Asymptotic independence}
 Let consider a process:
\begin{equation}
\left(S_{t},X_{t}\right)_{t=0}^{\infty}
\end{equation}
where 
\begin{itemize}
\item  $S_{t}$ is an unobservable hidden Markow chain  with $N$ states;
\item realizations of  process $X_{t}$ are  observed; 
\item  the  distribution of $X_{t}$ conditional  on history $\boldsymbol{{\cal R}}_{t-1}=\left(x_0,\ldots, x_{t-1}\right) $ has a form:
\begin{equation}
f\left (\boldsymbol{x}_{t}\mid\boldsymbol{{\cal R}}_{t-1};\boldsymbol{\theta}\right) 
 = \sum_{j=1}^{l}f\left (\boldsymbol{x}_{t}\mid S_{t}=j,\boldsymbol{{\cal R}}_{t-1};\boldsymbol{\theta}\right)P\left (S_{t}=j\mid\boldsymbol{{\cal R}}_{t-1};\boldsymbol{\theta}\right). \label{distr}\end{equation}
\end{itemize}
We prove the following lemma:
\begin{lemma}
Let consider a proces $\left(S_t, X_t\right)$ where the conditional distribution of $X_t\mid \boldsymbol{{\cal R}}_{t-1}$ is defined as \eqref{distr}. Let assume that    $S_t$ is an ergodic process. Then the random variables$\{X_t\}$ are asymptotically independent, i.e.
\begin{equation}\label{as-ind}
\lim_{s_1,\ldots,s_k \rightarrow \infty} \left|P(\bigcap_{\nu=1}^k\{X_{t+\sum_{\rho=0}^{\nu}s_{\rho}}\in A_{\rho}\})-\prod_{\nu=1}^{k}P(X_{t+\sum_{\rho=0}^{\nu}s_{\rho}}\in A_{\rho})\right|=0
\end{equation}
\end{lemma}

{\bf Proof:}
	At first we estimate the difference
\begin{multline*}\left|P(\bigcap_{\nu=1}^k\{X_{t+\sum_{\rho=0}^{\nu}s_{\rho}}\in A_{\rho}\})-\prod_{\nu=1}^{k}P(X_{t+\sum_{\rho=0}^{\nu}s_{\rho}}\in A_{\rho})\right|\\ \leq  \left|P(\bigcap_{\nu=1}^k\{X_{t+\sum_{\rho=0}^{\nu}s_{\rho}}\in A_{\rho}\})-P(\{X_{t+\sum_{\rho=0}^{k}s_{\rho}}\in A_{\rho}\})P(\bigcap_{\nu=1}^{k-1}\{X_{t+\sum_{\rho=0}^{\nu}s_{\rho}}\in A_{\rho}\})\right|\\+ \left|P(\{X_{t+\sum_{\rho=0}^{k}s_{\rho}}\in A_{\rho}\}P(\bigcap_{\nu=1}^{k-1}\{X_{t+\sum_{\rho=0}^{\nu}s_{\rho}}\in A_{\rho}\})-\prod_{\nu=1}^{k}P(X_{t+\sum_{\rho=0}^{\nu}s_{\rho}}\in A_{\rho})\right|=I_1+I_2\end{multline*}
\begin{multline*}I_1=\left|P(\bigcap_{\nu=1}^k\{X_{t+\sum_{\rho=0}^{\nu}s_{\rho}}\in A_{\rho}\})-P(\{X_{t+\sum_{\rho=0}^{k}s_{\rho}}\in A_{\rho}\})P(\bigcap_{\nu=1}^{k-1}\{X_{t+\sum_{\rho=0}^{\nu}s_{\rho}}\in A_{\rho}\})\right|\\=
\left|P(\bigcap_{\nu=1}^k\{X_{t+\sum_{\rho=0}^{\nu}s_{\rho}}\in A_{\rho}\}|\bigcap_{\nu=1}^{k-1}\{X_{t+\sum_{\rho=0}^{\nu}s_{\rho}}\in A_{\rho}\})P(\bigcap_{\nu=1}^{k-1}\{X_{t+\sum_{\rho=0}^{\nu}s_{\rho}}\in A_{\rho}\})\right.\\-\left.P(\{X_{t+\sum_{\rho=0}^{k}s_{\rho}}\in A_{\rho}\})P(\bigcap_{\nu=1}^{k-1}\{X_{t+\sum_{\rho=0}^{\nu}s_{\rho}}\in A_{\rho}\})\right|\\ \leq \left|P(\bigcap_{\nu=1}^k\{X_{t+\sum_{\rho=0}^{\nu}s_{\rho}}\in A_{\rho}\}\mid\{X_{t+\sum_{\rho=0}^{k-1}s_{\rho}}\in A_{\rho}\})-P(\{X_{t+\sum_{\rho=0}^{k}s_{\rho}}\in A_{\rho}\})\right|
\end{multline*}
Denote\begin{itemize}
\item$ T=t+\sum_{\rho=0}^{k-1}\tau{\rho}$
\item $\tau=\tau_k$
\item $A=A_k$
\item $B=A_{k-1}$
\end{itemize}

We notice, that \[P(X_{T+\tau}\in A, X_T\in B)=P(X_{T+\tau}\in A|X_T\in B)P(X_T\in B)\] 

Moreover for $A=\{s_1,\ldots, s_p\}\times \mathfrak{A}, B=\{s_0\}\times\mathfrak{B}$ \begin{equation}\label{war}
 P(X_{T+\tau}\in A |X_T\in B)=\sum_{\iota=1}^{p}p_{s_0,s_{\iota}}(\tau)\int_{\mathfrak{A}}f_{\iota}(x)dx
\end{equation}
and
\begin{multline} \label{bezwar}
 P(X_{T+\tau}\in A )=\sum_{\sigma \in\boldsymbol{{\cal S}}}P(X_{T+\tau}\in A \mid X_T\in \{\sigma\}\times \mathbb{R}^N)P(X_T\in \{\sigma\}\times \mathbb{R}^N) \\= \sum_{\sigma\in \boldsymbol{{\cal S}}}\sum_{\iota=1}^{p}p_{\sigma,s_{\iota}}(\tau)\int_{\mathfrak{A}}f_{\iota}(x)dxP(X_T\in \{\sigma\}\times \mathbb{R}^N)
\end{multline}
 From the ergodic theorem \cite{} follows, that \[\lim_{s\rightarrow \infty}p_{ij}(s)=p_j^*(s)\] and\[ |p_{ij}(s)-p_j^*(s)|\leq \alpha^s \] for some $\alpha <1$.
 
 From \eqref{bezwar}, \eqref{war} and the obvious equality \[\sum_{\sigma\in \boldsymbol{{\cal S}}}P(X_T\in \{\sigma\}\times \mathbb{R}^N)=1\] follows, that
 \begin{multline*}
 P(X_{T+\tau}\in A |X_T\in B)- P(X_{T+\tau}\in A ) \\ =\sum_{\iota=1}^{p}p_{s_0,s_{\iota}}(\tau)\int_{\mathfrak{A}}f_{\iota}(x)dx-\sum_{\sigma\in \boldsymbol{{\cal S}}}\sum_{\iota=1}^{p}p_{\sigma,s_{\iota}}(\tau)\int_{\mathfrak{A}}f_{\iota}(x)dxP(X_T\in \{\sigma\}\times \mathbb{R}^N) \\=\sum_{\iota=1}^{p}[p_{s_0,s_{\iota}}(\tau)-p_{s_{\iota}}^*(\tau)]\int_{\mathfrak{A}}f_{\iota}(x)dx +\sum_{\sigma\in \boldsymbol{{\cal S}}}\sum_{\iota=1}^{p}[p_{s_{\iota}}^*(\tau)-p_{\sigma,s_{\iota}}(\tau)]\int_{\mathfrak{A}}f_{\iota}(x)dxP(X_T\in \{\sigma\}\times \mathbb{R}^N)
 \end{multline*} and in consequence \[| P(X_{T+\tau}\in A |X_T\in B)- P(X_{T+\tau}\in A )|\leq 2 \alpha^{\tau}\]
 Hence $I_1\leq 2 \alpha^{\tau}$. Now. notice, that 
\begin{multline*}
 I_2=\left|P(\{X_{t+\sum_{\rho=0}^{k}\tau_{\rho}}\in A_{\rho}\})P(\bigcap_{\nu=1}^{k-1}\{X_{t+\sum_{\rho=0}^{\nu}\tau_{\rho}}\in A_{\rho}\})-\prod_{\nu=1}^{k}P(X_{t+\sum_{\rho=0}^{\nu}\tau_{\rho}}\in A_{\rho})\right|\\= P(\{X_{t+\sum_{\rho=0}^{k}\tau_{\rho}}\in A_{\rho}\})\left|P(\bigcap_{\nu=1}^{k-1}\{X_{t+\sum_{\rho=0}^{\nu}\tau_{\rho}}\in A_{\rho}\})-\prod_{\nu=1}^{k-1}P(X_{t+\sum_{\rho=0}^{\nu}\tau_{\rho}}\in A_{\rho})\right|\\ \leq \left|P(\bigcap_{\nu=1}^{k-1}\{X_{t+\sum_{\rho=0}^{\nu}\tau_{\rho}}\in A_{\rho}\})-\prod_{\nu=1}^{k-1}P(X_{t+\sum_{\rho=0}^{\nu}\tau_{\rho}}\in A_{\rho})\right|
\end{multline*}
 By simple induction we can conclude, that \[\left|P(\bigcap_{\nu=1}^k\{X_{t+\sum_{\rho=0}^{\nu}\tau_{\rho}}\in A_{\rho}\})-\prod_{\nu=1}^{k}P(X_{t+\sum_{\rho=0}^{\nu}\tau_{\rho}}\in A_{\rho})\right| \leq k\alpha^{\sum_{\rho=0}^{k}\tau_{\rho}}\] which completes the proof.
\section{Property of $\varepsilon$-independence}
Next, let  define the notion of $\varepsilon$-independence. This concept will be useful to prove central limit theorem. 

\begin{defn}
	The random variables sequence   $\{X_k\}_{k\in \mathbb{N}}$ are $\varepsilon$-independent, when for any  $n\in\mathbb{N}$ and any sets  $A_1,\ldots,A_n$ we have the following inequality:
	\begin{equation} \label{eps-indep}
	|P(X_1\in A_1, \ldots, X_n\in A_n)- P(x_1\in A_1)\ldots P(X_n\in A_n)|\leq\varepsilon. 
	\end{equation}
		\end{defn}
We  prove, that for $\varepsilon$-independent variables the following lemma is true. 
	\begin{lemma}
		When random variables  $X_1, \ldots, X_n$ are $\varepsilon$-independent, then  \[ |\varphi_{X_1+\ldots+X_n}(t)-\varphi_{X_1}(t)\ldots\varphi{X_n}(t)|\leq 2\varepsilon. \]
	\end{lemma}
	{\bf Proof:}
		From the formula \eqref{eps-indep} follows, that 
		\[|E(\mathbbm{1}_{A_1}(X_1)\ldots\mathbbm{1}_{A_n}(x_n))-E\mathbbm{1}_{A_1}(X_1)\ldots E\mathbbm{1}_{A_n}(x_n)|\leq \varepsilon.\]
		Since every coninuous function can be approximatad by simple functions consider at first  the real function  \begin{equation} \label{prosta}
		f=\sum_{j=1}^{m}c_{j}\mathbbm{1}_{A_{j}},
		\end{equation} where the sets $A_j$ are pairwise disjoint and $|c_j|\leq 1$.
In the first step, we estimate the real part of $\left(Ef(X_1)\cdots f(X_n)-Ef(X_1)\cdots Ef(X_n)\right)$. Since  $|c_{j}|\leq 1$ we get that:
\begin{multline*}
	\Re(Ef(X_1)\cdots f(X_n)-Ef(X_1)\cdots Ef(X_n)) \\ = \sum_{j_1,j_2,\ldots, j_n=1}^n\Re(c_{j_1}c_{j_2}\cdots c_{j_n})\left[P(X_1\in A_{j_1}, \ldots , X_n\in A_{j_n})-P(X_1\in A_{j_1})\cdots  P(X_n\in A_{j_n}) \right] \\ \leq  \sum_{j_1,j_2,\ldots, j_n=1}^n\left[P(X_1\in A_{j_1}, \ldots , X_n\in A_{j_n})-P(X_1\in A_{j_1})\cdots  P(X_n\in A_{j_n}) \right].  	
	 \end{multline*}
	 Since the sets $A_j$ are pairwise disjoint then:
	 \[\sum_{j_1,j_2,\ldots, j_n=1}^nP(X_1\in A_{j_1}, \ldots , X_n\in A_{j_n})=P(X_1\in A, \ldots, x_n\in A),\]
	 \[\sum_{j=1}^{m}P(X_k\in A_j)=P(X_k\in A),\] 
	 where \[ A=\bigcup_{j=1}^m A_j.\]
	 From it follows, that 
	 \[\Re(Ef(X_1)\cdots f(X_n)-Ef(X_1)\cdots Ef(X_n)) \leq P(X_1\in A, \ldots, x_n\in A)- P(X_1\in A) \cdots P(X_n\in A)\leq \varepsilon. \]
	 Analogously
 \[\Re(Ef(X_1)\cdots f(X_n)-Ef(X_1)\cdots Ef(X_n)) \geq -P(X_1\in A, \ldots, x_n\in A)+ P(X_1\in A) \cdots P(X_n\in A)\geq -\varepsilon \]
In the second  step, we estimate the imaginary  part of $\left(Ef(X_1)\cdots f(X_n)-Ef(X_1)\cdots Ef(X_n)\right)$. So, we have:
 \[\Im(Ef(X_1)\cdots f(X_n)-Ef(X_1)\cdots Ef(X_n)) \leq P(X_1\in A, \ldots, x_n\in A)- P(X_1\in A) \cdots P(X_n\in A)\leq \varepsilon \]
 \[\Im(Ef(X_1)\cdots f(X_n)-Ef(X_1)\cdots Ef(X_n)) \geq -P(X_1\in A, \ldots, x_n\in A)+ P(X_1\in A) \cdots P(X_n\in A)\geq -\varepsilon. \]
 In consequence we obtain:
  \[|Ef(X_1)\cdots f(X_n)-Ef(X_1)\cdots Ef(X_n)| \leq \sqrt{2}\varepsilon.  \]
    The last inequality is true for any $X_1, \ldots, X_n$  for any $n$ and for any $f$ Now, let fix  $t\in \mathbb{R}$, $n$ and $\eta>0$. Let define  
  the function $f_{\eta}$ which can be presented as \eqref{prosta} and satisfying the inequality \[|f_{\eta}(x)-e^{itx}|\leq \eta\] for every $x\in \mathbb{R}$.
  Let estimate the difference 
  \begin{equation}\varphi_{X_1+\ldots+X_n}(t)- Ef_{\eta}(X_1)\cdots f_{\eta}(X_n)
  \label{sacowanie}
  \end{equation}
  This difference is equal to \begin{multline}
  Ee^{itX_1}\cdots e^{itX_n}-Ef_{\eta}(X_1)\cdots f_{\eta}(X_n) \\ =  E(e^{itX_1}-f_{\eta}(X_1))e^{itX_2}\cdots e^{itX_n}+ Ef_{\eta}(X_1)(e^{itX_2}-f_{\eta}(X_2))e^{itX_3}\cdots e^{itX_n}+ \ldots \label{rownosc}
  \end{multline}
Each component of the  sum \eqref{rownosc} has the form $ EZ(e^{itX_k}- f_{\eta}(X_k)) $ where $|Z|\leq 1$. From this fact implies that  \[|Ee^{itX_1}\cdots e^{itX_n}-Ef_{\eta}(X_1)\cdots f_{\eta}(X_n)|\leq n\eta.\]
By conducting a similar reasoning, we obtain that:
 \[|Ee^{itX_1 +\cdots +itX_n}-E(f_{\eta}(X_1)\cdots f_{\eta}(X_n))|\leq n\eta.\]
 From it follows, that \begin{multline*}
 	  |Ee^{itX_1 +\cdots +itX_n}-Ee^{itX_1}\cdots e^{itX_n}| \\ \leq |Ee^{itX_1 +\cdots +itX_n}-E(f_{\eta}(X_1)\cdots f_{\eta}(X_n))|\\+Ef(X_1)\cdots f(X_n)-Ef(X_1)\cdots Ef(X_n)|+Ee^{itX_1}\cdots e^{itX_n}-Ef_{\eta}(X_1)\cdots f_{\eta}(X_n)| \\ \leq \sqrt{2}\varepsilon+ 2\eta n
 	  \end{multline*}

Choosing such $\eta$, that $2n\eta\leq 2\varepsilon - \sqrt{2}\varepsilon$ we complete the proof.

\section{Central Limit Theorem }	
Using the thesis of the above lemma, one can prove a central limit theorem for  $\varepsilon$-independent variables. 
\begin{theorem}
Let assume that for each  $ n $,    $  X_{{1n}}, X_{{2n}}, \ldots, X_{{r_{n}n}} $ are  $\varepsilon$-independent random variables  with expected value equal $ 0 $ and: 
\begin{equation}
\sum _{{k=1}}^{{r_{n}}}\mathbb{E}X_{{kn}}^{2}\xrightarrow{n\to\infty}1.\label{porownanie}	\end{equation}
Additionally, let us assume that Lindeberg's condition is fulfilled:
\[\sum _{{k=1}}^{{r_{n}}}\mathbb{E}X_{{kn}}^{2}1_{{\{|X_{{kn}}|>\eta\}}}\xrightarrow{n\to\infty}0\quad\mbox{for each }\eta>0.\]		
Let  $Y_n$ be a sum:
\[= X_{{1n}}+X_{{2n}}+\ldots+X_{{r_{n}n}}.\]
Then
\[\limsup|\varphi_{Y_n}(t)-e^{-\frac{1}{2}t^2}|<\varepsilon.\]
\end{theorem}

{\bf Proof:}
For Lindeberg Theorem (Loeve, 1977) implies that:
\[\lim_{n\rightarrow \infty}\varphi_{X_{1n}}(t)\ldots\varphi_{X_{r_nn}}(t)=e^{-\frac{1}{2}t^2}\]
From last Lemma implies that: \[|\varphi_{Y_n}(t)-\varphi_{X_{1n}}(t)\ldots\varphi_{X_{r_nn}}(t)|<\varepsilon.\]

Now, our goal is to prove the fact that asymptotically independent variables also  satisfy the Central Limit Theorem.
\begin{theorem}
	Assume, that
	\begin{enumerate}
		\item Random Variables $X_1, \ldots,X_n$ are asymptoticaly independent i.e., they  satisfy \eqref{as-ind}.
		\item $EX_i=0$ for $i=1,\ldots,n.$
			\item \begin{equation}
			E\left(\frac{1}{\sqrt{n}}\sum_{k=1}^n X_{i+k}\right)^2 \rightarrow 1  \text{   uniformly in }i. \label{wariancja}\end{equation}
		\end{enumerate}
		Then, the  sequence of the distributions of the random variables \[ \frac{1}{\sqrt{v}}(X_1+\cdots+X_n)\]
		 tend to the standard normal distribution.
\end{theorem}

	{\bf Proof:}
		At first we fix some arbitrary $ \varepsilon >0 $. From \eqref{as-ind} follows, that  there exists such $m$, that 	\[|P(X_{k_1}\in A_1,\ldots, X_{k_r}\in A_r)=P(X_{k_1}\in A_1)\ldots P(x_{k_r}\in A_r)|\leq \varepsilon,\] where $k_{j+1}>k_j+m.$
		\\
		Fix $n$ and let $0<\alpha<\frac{1}{4}$. Let define  $ k=[n^{\alpha}] $ and $ \nu=\left[\frac{k}{n}\right] $, so clearly $k\leq n^{\alpha}$ and  $n=k\nu+r$. \\
		Let define:
		\[U_i=X_{ik-k+1}+\ldots + X_{ik-m}.\]
		Because   $(i+1)k-k+1-(ik-m)=m+1$ so $U_i$ are $\varepsilon$-independent.  Let us consider a sum:
		 \[X_1+\ldots+X_n. \]
		 This sum we can seperate  into two parts: a sum of  $U_1+\ldots, U_{\nu}$ and the rest. We can notice that each component  $U_i$ includes  $k+1-m$  elements, so  the sum  $U_1+\ldots +U_{\nu}$ consists of  $(k+1-m)\nu= n-r-(m-1)\nu$ elements. So the rest includes  $r+(m-1)\nu$ ingredients,  which we denote as $Z_1^{(n)},\ldots,Z_p^{(n)}$, for fixed  $n$. So:
		  \[X_1+\ldots +X_n=U_1+\ldots+ U_{\nu}  +Z_1^{(n)}+\ldots Z_p^{(n)}=\sqrt{n}\mathfrak{U}_n+\sqrt{n}\mathfrak{Z}_n.\]
		 
		 From  Schwarz and  H\"older inequality (Vuong, 1989) implies that:
					 \[E(Z_1^{(n)}+\ldots +Z_p^{(n)})^2\leq p^2R^2,\]
			where $R=E|X_i|^3$, so: 
			 \[E(\frac{1}{\sqrt{n}}(Z_1^{(n)}+\ldots+Z_p^{(n)}))^2\leq \frac{p^2}{n}R^2.\]
			From a fact that  $p=r+(m-1)\nu\leq m\nu \leq m\frac{k}{n}\leq m\frac{n^{\alpha}}{n}$ implies that 
			  $\frac{p^2}{n}\leq m^2\frac{n^{2\alpha}}{n^3}$. 
			  So, the sum  \[\sum_{n=1}^{\infty}\frac{p^2}{n}R^2\] 
			  is consistent, so with probability $1$, we have:
			  \begin{equation} \label{z}
			  \lim_{n\rightarrow \infty}\frac{1}{\sqrt{n}}(Z_1^{(n)}+\ldots+Z_p^{(n)})=0.   
			  \end{equation}
			Now, we need estimate the 
			 \[\frac{1}{\sqrt{\nu}}(U_1^{(n)}+\ldots+U_{\nu}^{(n)})\] \[\frac{1}{\sqrt{n}}(U_1^{(n)}+\ldots+U_{\nu}^{(n)})=\frac{\sqrt{k\nu}}{\sqrt{n}}\frac{1}{\sqrt{k\nu}}(U_1^{(n)}+\ldots+U_{\nu}^{(n)})=\frac{\sqrt{k\nu}}{\sqrt{n}}\frac{1}{\sqrt{\nu}}(\frac{1}{\sqrt{k}}U_1^{(n)}+\ldots+\frac{1}{\sqrt{k}}U_{\nu}^{(n)}).
			 \] Moreover
			\[U_i^{(n)}=X_{ik-k+1}+\ldots + X_{ik-m},\] from which follows, that
			\[E\left(\frac{1}{\sqrt{k}}U_i^{(n)}\right)^2=\frac{k-m}{k}E\left(\frac{1}{\sqrt{k-m}}U_i^{(n)}\right)^2.\]
			From the assumption \eqref{wariancja}follows, that
			\[\lim_{k\rightarrow \infty}E\left(\frac{1}{\sqrt{k}}U_i^{(n)}\right)^2=1.\] In consequencewe obtain that:
						\[\lim_{n\rightarrow \infty}\frac{1}{\nu}\left(\sum_{i=1}^{\nu}E\left(\frac{1}{\sqrt{k}}U_i^{(n)}\right)^2\right)=1.\]
	Having above estimates we go to the next step of the proof. We have to estimate the difference:
	\begin{equation}\left| \varphi_{\frac{1}{\sqrt{n}}(X_1+\ldots+X_n)}(t)- e^{-\frac{t^2}{2}}\right|.
	  \end{equation}
So, 
	\begin{multline*}\left| \varphi_{\frac{1}{\sqrt{n}}(X_1+\ldots+X_n)}(t)- e^{-\frac{t^2}{2}}\right|=
	\left| \varphi_{\mathfrak{U}_n+\mathfrak{Z}_n}(t)- e^{-\frac{t^2}{2}}\right|=\\
	 = \left| \varphi_{\mathfrak{U}_n+\mathfrak{Z}_n}(t) -\varphi_{\mathfrak{U}_n}(t)+\varphi_{\mathfrak{U}_n}(t) - e^{-\frac{t^2}{2}}\right| \leq 	
	 	 \left|\varphi_{\mathfrak{U}_n+\mathfrak{Z}_n}(t) -\varphi_{\mathfrak{U}_n}(t)\right|+\left|\varphi_{\mathfrak{U}_n}(t) - e^{-\frac{t^2}{2}}\right|\\
	 	  \leq E\left|e^{it(\mathfrak{U}_n+\mathfrak{Z}_n)} -e^{it\mathfrak{U}_n}\right|+\left|\varphi_{\mathfrak{U}_n}(t) - e^{-\frac{t^2}{2}}\right|=E\left|e^{it(\mathfrak{U}_n}(e^{it\mathfrak{Z}_n)} -1)\right|+\left|\varphi_{\mathfrak{U}_n}(t) - e^{-\frac{t^2}{2}}\right|\\ \leq E\left|e^{it\mathfrak{Z}_n} -1\right|+\left|\varphi_{\mathfrak{U}_n}(t) - e^{-\frac{t^2}{2}}\right|
	  \end{multline*}
	  From \eqref{z} follows, that \[\lim_{n \rightarrow \infty }E\left|e^{it\mathfrak{Z}_n} -1\right|=0\] and from lemma follows, that \[\limsup_{n \rightarrow \infty }\left|\varphi_{\mathfrak{U}_n}(t) - e^{-\frac{t^2}{2}}\right|\leq 2 \varepsilon\]
	  Then
	  \[\limsup_{n \rightarrow \infty }\left| \varphi_{\frac{1}{\sqrt{n}}(X_1+\ldots+X_n)}(t)- e^{-\frac{t^2}{2}}\right|\leq 2 \varepsilon\]
	  Since the last estimation can be proved for any $\varepsilon$ 
	   \[\lim_{n \rightarrow \infty }\left| \varphi_{\frac{1}{\sqrt{n}}(X_1+\ldots+X_n)}(t)- e^{-\frac{t^2}{2}}\right|=0\]
	   which completes the proof.

\section*{Conclusion}	
I this paper we prove a Central Limit Theorem. The assumptions to this thorem is not restrictive and they are not difficult to apply in the  practical.  We proved that when  a Markov chain is ergodic and random variables fullfiled Lindeberg's condition  then the Central Limit Theorem is true.	
\section*{References}
\[\]
 Garcia, R., (1998), Asymptotic null distribution of the likelihood ratio test in 
 
  Markov switching
models. International Economic Review 39, 763-88. \\
 Hamilton, J., (1989), A new approach to the economic analysis of nonstationary 
 
 time series
and the business cycle. Econometrica 57, 357-384.\\
  Hamilton, J., (1996), Specification testing in Markov-switching time-series models.
  
   Journal
of Econometrics 70, 127-157. \\
Hansen, B. E., (1992), The likelihood ratio test under non-standard conditions.

 Journal of Applied Econometrics 7, S61-82.\\ 
 Hoeffding W., Robbins H. ,(1994), The Central Limit Theorem for Dependent 
 
 Random Variables. In: Fisher N.I., Sen P.K. (eds) The Collected Works 
 
 of Wassily Hoeffding. Springer Series in Statistics. Springer, New York. \\
 Kim, C.-J., (2009), Markov-switching models with endogenous explanatory 
 
 variables II:  A two-step MLE procedure. Journal of Econo
 
 metrics 148, 46-55. \\
 Kim, C.-J., Nelson, C., (1999),  State-Space Models with Regime Switching. 
   MIT Press, Cambridge, MA.\\
  Loeve, M.,(1977), Probability theory, Springer. \\
  Timmermann, A., (2000), Moments of Markov switching models. Journal of Econo  
  
  metrics 96, 75-111. \\
   Vuong Q. H., (1989), 
   Likelihood ratio tests for model selection and non-nested 
   
   hypotheses, Econometrica 57, 307 - 333.

\end{document}